\documentclass[a4paper]{article}     % `article' with A4 paper size option
\usepackage{amsthm}                 % theorem extensions
\usepackage{amssymb}
\usepackage{mathrsfs}
\usepackage{dsfont}

\usepackage{amsmath}
\usepackage{ae}
\usepackage[latin1]{inputenc}
\usepackage[T1]{fontenc}
\usepackage[all]{xy}

\newtheorem{theorem}{Theorem}
\newtheorem*{theorem*}{Theorem}
\newtheorem{corollary}{Corollary}
\newtheorem{proposition}{Proposition}
\newtheorem{lemma}{Lemma}

\newcommand{\ind}{\textnormal{ind}}

\newcommand{\res}{\textnormal{res}}
\newcommand{\im}{\textnormal{im}}

\newcommand{\mo}{\textnormal{mod}}

\newcommand{\ke}{\textnormal{ker}}

\newcommand{\ann}{\textnormal{ann}}

\newcommand{\tot}{\textnormal{tot}}
\newcommand{\coke}{\textnormal{coker}}
\newcommand{\AAA}{{\mathscr{A}}}
\newcommand{\BB}{{\mathscr{B}}}
\newcommand{\CC}{{\mathds{C}}}

\newcommand{\RR}{{\mathds{R}}}

\newcommand{\NN}{{\mathds{N}}}

\newcommand{\calB}{{\mathcal{B}}}

\newcommand{\calG}{{\cal G}}

\newcommand{\OCn}{\mathscr{O}}

\newcommand{\ICn}{\ind_{\CC^n,0}}
\newcommand{\IVn}{\ind_{V,0}}

\newcommand{\RCn}{\res_{\CC^n,0}}

\newcommand{\omv}{\Omega_{V,0}}
\newcommand{\ov}{\mathscr{O}_{V,0}}
\newcommand{\omc}{\Omega}

\newcommand{\ACn}{\mathscr{A}}
\newcommand{\BCn}{\mathscr{B}}
\newcommand{\BCnp}{\mathscr{B}'}

\begin{document}

\title{On the index of a vector field tangent to a hypersurface with non-isolated zero in the embedding space}
\author{Oliver Klehn\\
\parbox{9cm}{\small
 \begin{center} Institut f\"{u}r Mathematik, Universit\"{a}t Hannover, \\
        Postfach 6009, D-30060 Hannover, Germany \\
        E-mail: klehn@math.uni-hannover.de
\end{center}}
\date{}
}

\maketitle

\begin{abstract}
We give a generalization of an algebraic formula of Gomez-Mont for the index of a vector field with isolated zero in
$(\CC^n,0)$ and tangent to an isolated hypersurface singularity. We only assume that the vector field has an isolated zero on the
singularity here.
\end{abstract}

\section{Introduction}
\label{intro}
Let $(V,0)=(\{ f=0\} ,0)\subseteq (\CC^n,0)$ be an isolated singularity of a hypersurface and
\[
X=X_1\frac{\partial}{\partial z_1}+\dots +X_n\frac{\partial}{\partial z_n}
\]
the germ of a holomorphic vector
field on $(\CC^n,0)$ tangent to $(V,0)$, say $Xf=cf$, with an isolated zero on $(V,0)$. In \cite{gsv} is defined a notion of an index $\IVn X$ of
this vector field, called GSV index and generalizing the Poincar\'{e}-Hopf index. By contracting germs of Kähler forms on
$(V,0)$ with $X$ one obtains two complexes
\[
(\omv^*,X)\colon 0\to\omv^n\stackrel{X}{\to}\cdots\stackrel{X}{\to}\omv^1\stackrel{X}{\to}\ov\to 0
\]
\[
(\omv,X)\colon 0\to\omv^{n-1}\stackrel{X}{\to}\cdots\stackrel{X}{\to}\omv^1\stackrel{X}{\to}\ov\to 0
\]
which have both finite dimensional homology spaces, see \cite{g}. We set
\begin{eqnarray*}
h_i^*:=\dim_{\CC}H_i(\omv^*,X) && h_i:=\dim_{\CC}H_i(\omv,X).
\end{eqnarray*}
It has been proved in \cite{g} that the topologically defined GSV index is given as
\[
\IVn X=\chi(\omv ,X)=\sum_{i=0}^{n-1}(-1)^ih_i .
\]
Further the numbers $h_i$ and $h_i^*$ are computed under the condition that $X$ has an isolated zero on $(\CC^n,0)$. We
want to extend these computations to the case that $X$ has an isolated zero only on $(V,0)$ using the same methods as in \cite{g}.

Note that one can assume after a linear change of coordinates in $(\CC^n,0)$ that $(X_1,\dots ,X_{n-1},f)$ is a regular
$\OCn$-sequence. Throughout this paper we always assume the coordinates to be chosen in this way. We denote by
\[
\ACn :=\frac{\OCn}{(\partial f/\partial z_1,\dots ,\partial f/\partial z_n)}
\]
the Milnor algebra of $f$ at $0$ and we set
\begin{eqnarray*}
\BCn :=\frac{\OCn}{(X_1,\dots ,X_n)} && \BCnp :=\frac{\OCn}{(X_1,\dots ,X_{n-1})}.
\end{eqnarray*}
We simply write $\OCn$ resp. $\omc^i$ instead of $\mathscr{O}_{\CC^n,0}$ resp. $\Omega^i_{\CC^n,0}$.
$\BCnp$ is infinite dimensional and $\BCn$ not necessarily finite dimensional. Our main result is the following:

\begin{theorem}
\label{alg2}
The $\CC$-dimensions of the homology spaces of $(\omv ,X)$ are given by
\[h_0^*=\dim_{\CC}\frac{\BCn}{(f)}, \]
\[h_1^*=\dim_{\CC}\frac{\ann_{\BCn}(f)}{(c)}+\dim_{\CC}\frac{\ann_{\BCnp}(X_n)}{\ann_{\BCnp}(X_n)\cap\BCnp
(f,\partial f/\partial z_n)}, \]
\[h_2^*=\dots =h_n^*=\dim_{\CC}\frac{\ann_{\ACn}(f)}{(c)}=\dim_{\CC}\frac{\ann_{\ACn}(c)}{(f)}=:\lambda . \]
\end{theorem}

\section{The homological index}
\label{homolog}
The first steps are easy generalizations of the computations in \cite{g}. Consider the complex
\[
(\omc^*,X)\colon 0\to\omc^n\stackrel{X}{\to}\cdots\stackrel{X}{\to}\omc^1\stackrel{X}{\to}\OCn\to 0.
\]

\begin{lemma}
\label{kos1}
$H_j(\omc^*,X)=0 \text{ for } j\geq 2$\\
$H_1(\omc^*,X)\cong \ann_{\BCnp}(X_n)$, ($h_1dz_1+\dots +h_ndz_n\to h_n$)\\
$H_0(\omc^*,X)\cong\BCn$.
\end{lemma}

\begin{proof}
It follows immediately from the dual version of Corollary 17.12 of \cite{e}.
\end{proof}

Multiplication with $c$ gives us a map of complexes
\[
\cdot c\colon (\omc^*,X)\to (\omc^*,X).
\]
We denote by $\calB$ the mapping cone associated to this map, see \cite{e}. We have
\[
\calB\colon
0 \leftarrow \OCn\oplus 0\stackrel{\beta}{\leftarrow} \omc^1\oplus\OCn\stackrel{\beta}{\leftarrow}
\dots\stackrel{\beta}{\leftarrow} \omc^n\oplus\omc^{n-1}\stackrel{\beta}{\leftarrow}
0\oplus\omc^n\stackrel{\beta}{\leftarrow} 0
\]
where
\[
\beta (\omega ,\eta ):=
\begin{pmatrix}
X & (-1)^jc\\
0 & X\\
\end{pmatrix}
\begin{pmatrix}
\omega \\ \eta \\
\end{pmatrix}
:=
\begin{pmatrix}
X(\omega )+ (-1)^jc\eta\\ X(\eta )\\
\end{pmatrix},
\]
and $j$ is the degree of $(\omega ,\eta )\in\omc^j\oplus\omc^{j-1}$.

\begin{lemma}
\label{kos2}
$H_k(\calB )=0$ for $k\geq 3$\\
$H_2(\calB )\cong\ann_{\BCnp}(X_n)\cap\ann_{\BCnp}(c)$\\
$H_1(\calB )\cong\frac{\ann_{\BCnp}(X_n)}{c\cdot\ann_{\BCnp}(X_n)}\oplus\ann_{\BCn}(c)$\\
$H_0(\calB )\cong\frac{\BCn}{(c)}$.
\end{lemma}
\begin{proof}
$H_0(\calB )\cong\frac{\BCn}{(c)}$ is obvious. We have an exact sequence of complexes
\[
\xymatrix{
  & 0\ar[d]    &    0\ar[d]    &                 &          0\ar[d]         &          &       \\
0 &\OCn\ar[l]\ar[d] & \omc^1\ar[l]_X\ar[d] & \dots \ar[l]_X & \omc^n\ar[l]_X\ar[d] &
0\ar[l]\ar[d]\\
0 &\OCn\oplus 0\ar[l]\ar[d] & \omc^1\oplus\OCn\ar[l]_{\beta}\ar[d]
& \dots \ar[l]_{\beta} & \omc^n\oplus\omc^{n-1}\ar[l]_{\beta}\ar[d] &
0\oplus\omc^n\ar[l]_{\beta}\ar[d] & 0\ar[l]\\
&0 &\OCn\ar[l]\ar[d] & \dots\ar[l]_X & \omc^{n-1} \ar[l]_X\ar[d]
& \omc^n\ar[l]_X\ar[d] &
0\ar[l]\\
&  & 0 && 0 & 0 & \\
}
\]
From the long exact sequence of homology we get $H_k(\calB )=0$ for $k\geq 3$ and the non-trivial part of the sequence is
\[
0\xrightarrow{}H_2(\calB )\xrightarrow{p_2}H_1(\Omega^*,X)\xrightarrow{\cdot c}H_1(\Omega^*,X)\xrightarrow{i_1}H_1(\calB
)\xrightarrow{p_2}\BCn\xrightarrow{\cdot c}\BCn\xrightarrow{}\frac{\BCn}{(c)}\xrightarrow{} 0,
\]
where $i_1$ is the inclusion to the first factor and $p_2$ the projection onto the second factor.
Now one has
\[
H_2(\calB )\cong\im p_2=\ker c\cong\ann_{\BCnp}(X_n)\cap\ann_{\BCnp}(c)
\]
by Lemma \ref{kos1} and the isomorphism given there. Further we get
\[
H_1(\calB )\cong\ker p_2\oplus\im p_2.
\]
We have $\im p_2=\ker c=\ann_{\BCn}(c)$ and
\[
\ker p_2=\im i_1\cong\frac{H_1(\Omega^*,X)}{\ker i_1}=\frac{H_1(\Omega^*,X)}{\im c}\cong
\frac{\ann_{\BCnp}(X_n)}{c\cdot\ann_{\BCnp}(X_n)}.
\]
\end{proof}

\begin{lemma}
\label{kos3}
\[
H_1(\omv^*,X)\cong\frac{\ann_{\BCn}(f)}{(c)}\oplus\frac{\ann_{\BCnp}(X_n)}{\ann_{\BCnp}(X_n)\cap\BCnp (f,\partial f/\partial z_n)}
\]
\end{lemma}

\begin{proof}
We have an exact sequence of complexes
\[
\xymatrix{
& 0\ar[d]  &  0\ar[d] &   &   0\ar[d]  &  \\
0 & \mathscr{R}^0\ar[l] \ar[d]& \mathscr{R}^1\ar[l]_X\ar[d]  &  \dots\ar[l]_X & \mathscr{R}^n\ar[l]_X \ar[d]& 0\ar[l] \\
0 & \OCn\ar[l] \ar[d]& \omc^1\ar[l]_X\ar[d]  &  \dots\ar[l]_X & \omc^n\ar[l]_X \ar[d]& 0\ar[l] \\
0 & \ov\ar[l] \ar[d]& \omv^1\ar[l]_X\ar[d]  &  \dots\ar[l]_X & \omv^n\ar[l]_X \ar[d]& 0\ar[l] \\
& 0 &  0&   &   0 &  \\
}
\]
where $\mathscr{R}^k:=f\omc^k +df\wedge\omc^{k-1}$.
From Lemma \ref{kos1} one gets the last part of the long exact sequence in homology
\begin{multline*}
0\xrightarrow{} H_2(\omv^*,X)\xrightarrow{\delta_2}H_1(\mathscr{R},X )\xrightarrow{i_1}H_1(\Omega^*,X)\xrightarrow{p_1}
H_1(\omv^*,X)\xrightarrow{\delta_1}\\ \xrightarrow{\delta_1}H_0(\mathscr{R},X )\xrightarrow{i_0}\BCn\xrightarrow{p_0}\frac{\BCn}{(f)}\xrightarrow{}0,
\end{multline*}
where $\delta_1$ and $\delta_2$ are the connecting homomorphisms. This gives us
\[
H_1(\omv^*,X)\cong\ker\delta_1\oplus\im\delta_1.
\]
One has $\im \delta_1 =\ker i_0$, $\mathscr{R}^0=f\OCn$ and $\mathscr{R}^1=f\Omega^1+df\OCn$. For
\[
\omega =
f(h_1dz_1+\dots +h_ndz_n)+hdf
\]
we get
\[
X(\omega )=f(X_1h_1+\dots +X_nh_n)+hcf
\]
and therefore
\[
\ker i_0=\frac{\{ gf\in\OCn :gf\in I_n \}}{fI_n+\OCn (cf)}.
\]
We have set $I_n:=\OCn (X_1,\dots ,X_n)$ here.
On the other hand multiplication gives an isomorphism
\[
\cdot f\colon \frac{\ann_{\BCn}(f)}{(c)}\to\frac{\{ gf\in\OCn :gf\in I_n \}}{fI_n+\OCn (cf)}
\]
and this proves
\[
\im\delta_1\cong\frac{\ann_{\BCn}(f)}{(c)}.
\]
Consider $\ker\delta_1$ now. Then we have
\[
\ker\delta_1=\im p_1\cong\frac{H_1(\Omega^*,X)}{\ker p_1}=\frac{H_1(\Omega^*,X)}{\im i_1}.
\]
We have an isomorphism
\[
\phi\colon H_1(\Omega^*,X)\to \ann_{\BCnp}(X_n)
\]
given as $\phi (h_1dz_1+\dots +h_ndz_n)=h_n$ with $X_1h_1+\dots +X_nh_n=0$. To prove
\[
\ker\delta_1\cong\frac{\ann_{\BCnp}(X_n)}{\ann_{\BCnp}(X_n)\cap\BCnp (f,\partial f/\partial z_n)}
\]
we have to show
\[
\phi\circ i_1(H_1(\mathscr{R},X ))=\ann_{\BCnp}(X_n)\cap\BCnp (f,\partial f/\partial z_n).
\]
" $\subseteq$ ": Let $\omega\in H_1(\mathscr{R},X )$, say $\omega =f(h_1dz_1+\dots +h_ndz_n)+hdf$ with $X(\omega )=0$ and so we have
$f(X_1h_1+\dots +X_nh_n) + hcf=0$. $\omega$ is mapped to $fh_n+h\partial f/\partial z_n$. Modulo $I_{n-1}$ we have
\[
(fh_n+h\partial f/\partial z_n)X_n=fh_nX_n+hcf=0
\]
and this shows the inclusion. Here we have set $I_{n-1}:=\OCn (X_1,\dots ,X_{n-1})$.\\
"$\supseteq $": Let $h\in\ann_{\BCnp}(X_n)\cap\BCnp (f,\partial f/\partial z_n)$ be a representative, say
\[
hX_n=g_1X_1+\dots +g_{n-1}X_{n-1}
\]
and $h=s_1f+s_2\partial f/\partial z_n$ with $s_1,s_2\in\OCn$. It follows that
\begin{equation*}
\begin{split}
s_1fX_n+s_2\partial f/\partial z_n X_n &=f(s_1X_n+s_2c)-s_2(\partial f/\partial z_1 X_1+\dots +\partial f/\partial z_{n-1} X_{n-1})\\
&=g_1X_1+\dots +g_{n-1}X_{n-1}
\end{split}
\end{equation*}
and therefore $f(s_1X_n+s_2c)$ is contained in $I_{n-1}$. This means that $s_1X_n+s_2c$ is contained in $I_{n-1}$, because $f$ is not a zero
divisor in $\BCnp$. Let
\[
s_1X_n+s_2c=l_1X_1+\dots +l_{n-1}X_{n-1}.
\]
Then one has
\[
g_1X_1+\dots +g_{n-1}X_{n-1}=fl_1X_1+\dots +fl_{n-1}X_{n-1}-s_2(\partial f/\partial z_1 X_1+\dots +\partial f/\partial z_{n-1} X_{n-1}).
\]
Set
\[
\omega :=-(fl_1-s_2\partial f/\partial z_1)dz_1-\dots -(fl_{n-1}-s_2\partial f/\partial z_{n-1})dz_{n-1}+hdz_n.
\]
It follows
\[
\omega =f(-l_1dz_1-\dots -l_{n-1}dz_{n-1}+s_1dz_n)+s_2\cdot df\in\mathscr{R}^1
\]
and
\begin{equation*}
\begin{split}
X(\omega ) & =f(-X_1l_1-\dots -X_{n-1}l_{n-1}+s_1X_n)+s_2cf\\
& =f(-s_1X_n-s_2c+s_1X_n)+s_2cf\\
& =0\\
\end{split}
\end{equation*}
Therefore $w\in H_1(\mathscr{R},X )$ with $\phi\circ i_1(\omega )=h$.
\end{proof}

\subsection{Proof of Theorem \ref{alg2}}
\paragraph{The double complex $\mathcal{G}$}

We use the construction as in \cite{g}. Consider the double complex $\calG$:
\[
\xymatrix{
&&&&& 0\ar[d] & \\
&&&& 0\ar[d] & \OCn\oplus 0\ar[d]^{\alpha}\ar[l] & \dots\ar[l]_{\beta}\\
&&& \dots &\OCn\oplus 0\ar[d]^{\alpha}\ar[l] & \omc^1\oplus\OCn\ar[l]_{\beta}\ar[d]^{\alpha} &\dots\ar[l]_{\beta}\\
&&& \dots &\omc^1\oplus\OCn\ar[d]^{\alpha}\ar[l] & \omc^2\oplus\omc^1 \ar[l]_{\beta}\ar[d]^{\alpha} &\dots\ar[l]_{\beta}\\
&& 0\ar[d] & \dots\ar[l] & \dots\ar[d]^{\alpha}  & \dots\ar[d]^{\alpha} & \dots \\
& 0\ar[d] & \OCn\oplus 0\ar[l]\ar[d]^{\alpha} & \dots\ar[l]_{\beta}& \omc^{n-1}\oplus\omc^{n-2}\ar[d]^{\alpha}\ar[l] & \omc^n\oplus\omc^{n-1}\ar[l]_{\beta}\ar[d]^{\alpha} & \dots\ar[l]_{\beta}\\
0 & \OCn\oplus 0\ar[l]\ar[d]^f & \omc^1\oplus\OCn\ar[l]_{\beta}\ar[d]^{\bar{\alpha}}& \dots\ar[l]_{\beta} & \omc^n\oplus\omc^{n-1}\ar[d]^{\bar{\alpha}}\ar[l] & 0\oplus\omc^n\ar[d]\ar[l]_{\beta} & 0\ar[l]_{\beta}\\
0& \OCn\ar[d]\ar[l] & \omc^1\ar[d]\ar[l]_X & \dots\ar[l]_X & \omc^n\ar[d]\ar[l] & 0\ar[l] & \\
& 0& 0 && 0&&\\
}
\]
Numeration: $\calG_{i,0}=\Omega^i$, $\calG_{i,1}=\Omega^i\oplus\Omega^{i-1}=\calG_{i+k,1+k}$ for $k\geq 1$ and
$\calG_{i,j}=0$ for $i<0$ or $j<0$. The mappings are $\bar{\alpha}(\omega ,\eta ):=\pm f\omega+df\wedge\eta$,
\[
\alpha (\omega ,\eta ):=\begin{pmatrix}
df & 0\\ \pm f & df
\end{pmatrix}\begin{pmatrix} \omega \\ \eta \end{pmatrix}:=
\begin{pmatrix}
df\wedge\omega \\ \pm f\omega +df\wedge\eta
\end{pmatrix}
\]
and
\[
\beta (\omega ,\eta ):=\begin{pmatrix}
X & \pm c\\ 0 & X
\end{pmatrix}\begin{pmatrix} \omega \\ \eta \end{pmatrix}:=
\begin{pmatrix}
X( \omega )\pm c\eta\\ X(\eta )
\end{pmatrix}.
\]
In \cite{g} it is shown that
\[
H_i(\tot\calG )\cong H_i(\omv^*,X) \textit{ for } 0\leq i\leq n
\]
and
\[
\dim_{\CC}H_i(\tot\calG )=\lambda \textit{ for } i>n,
\]
where $\tot\calG$ denotes the total complex of $\calG$ and $\lambda$ is defined as in Theorem \ref{alg2}.
The argument there does not use at this stage that $(X_1,\dots ,X_n)$ is a regular sequence.
The homology spaces of the total complex are also obtained by analizing the spectral sequence
$(E,d)$ of $\mathcal{G}$ obtained by taking first the $\alpha$-homology and then the $\beta$-homology.

\paragraph{The spectral sequence $(E,d)$}
Let us analyze $(E,d)$. We get the following
$^2E$-terms by Lemma \ref{kos1} and Lemma~\ref{kos2}:
\begin{align*}
^2E_{0,0}&\cong\coke \left( \frac{\BB}{(c)}\xrightarrow{\cdot f}\BB\right) \cong\frac{\BB}{(f)}\\
^2E_{0,1}&\cong\ke \left( \frac{\BB}{(c)}\xrightarrow{\cdot f}\BB\right)  \cong\frac{\ann_{\BB}(f)}{(c)}\\
^2E_{1,0}&\cong\coke \left( H_1(\calB)\xrightarrow{\bar{\alpha}}H_1(\Omega^*,X)\right) \\
^2E_{1,1}&\cong\text{ homology }\left( \frac{\BB}{(c)}\xrightarrow{\alpha}H_1(\calB)\xrightarrow{\bar{\alpha}}H_1(\Omega^*,X)\right) \\
^2E_{l+1,l}&\cong\coke \left( H_1(\calB )\xrightarrow{\alpha}H_2(\calB )\right) \\
^2E_{l,l+1}&\cong\ke \left( \frac{\BB}{(c)}\xrightarrow{\alpha}H_1(\calB )\right) \\
^2E_{l+1,l+1}&\cong\text{ homology }\left( \frac{\BB}{(c)}\xrightarrow{\alpha}H_1(\calB)\xrightarrow{\alpha}H_2(\calB)\right) \\
\end{align*}
where $l\geq 1$ and all other $^2E$-terms are trivial. The differential on $^2E$ is now given as
\[
d_2\colon ^2E_{q-1,p+2}\to ^2E_{q,p}.
\]
Therefore the spectral sequence degenerates at the second level and we are able to compute the total homology of the double complex. First note that
\[
\ke \left( H_1(\calB )\xrightarrow{\bar{\alpha }}H_1(\Omega^*,X)\right) =\ke \left( H_1(\calB )\xrightarrow{\alpha}H_2(\calB )\right)
\]
by the isomorphisms given in Lemma \ref{kos1} and Lemma \ref{kos2}. Therefore we get
\[
\dim_{\CC}H_2(\tot\calG )=\dim_{\CC }H_i(\tot\calG )=\lambda
\]
for $i\geq 2$ even. Similarly we get
\[
\dim_{\CC}H_i(\tot\calG )=\lambda
\]
for $i\geq 3$ odd. Comparing this with the result on $H_i(\tot\mathcal{G})$ in \cite{g} as stated above and
using Lemma \ref{kos3} for $h_1^*$ Theorem \ref{alg2} follows.

\begin{corollary}
\label{alg}
The following algebras are finite dimensional and for $n$ even
\begin{multline*}
\IVn (X)=\dim_{\CC}\frac{\BB}{(f)}-\dim_{\CC}\frac{\ann_{\BB}(f)}{(c)}-
\dim_{\CC}\frac{\ann_{\BB'}(X_n)}{\ann_{\BB'}(X_n)\cap\BB' (f, \partial f/\partial z_n)}\\+\dim_{\CC}\frac{\AAA}{(c)}-\dim_{\CC}\AAA
\end{multline*}
and for $n$ odd
\[
\IVn (X)=\dim_{\CC}\frac{\BB}{(f)}-\dim_{\CC}\frac{\ann_{\BB}(f)}{(c)}-
\dim_{\CC}\frac{\ann_{\BB'}(X_n)}{\ann_{\BB'}(X_n)\cap\BB' (f, \partial f/\partial z_n)}+\dim_{\CC}\frac{\AAA}{(f)}.
\]
\end{corollary}

\begin{proof}
First note that $\omv^n\cong\frac{\AAA}{(f)}$ and so we get
\begin{equation*}
\begin{split}
h_{n-1} &=h_{n-1}^*+\dim_{\CC}X(\omv^n)\\
            &=h_{n-1}^*+\dim_{\CC}\frac{\AAA}{(f)}-h_n^*.
\end{split}
\end{equation*}
Now it follows
\[
\sum_{\nu =0}^{n-1}(-1)^{\nu}h_i=\sum_{\nu =0}^n(-1)^{\nu}h_i^*+(-1)^{n-1}\dim_{\CC}\frac{\AAA}{(f)}
\]
and we get the claim for $n$ odd and for $n$ even
\[
\IVn (X)=h_0^*-h_1^*-\dim_{\CC}\frac{\AAA}{(f)}+\dim_{\CC}\frac{\ann_{\AAA}(f)}{(c)}.
\]
Since we have an exact sequence
\[
0\xrightarrow{} \ann_{\AAA}(f)\xrightarrow{} \AAA \xrightarrow{\cdot f} \AAA \xrightarrow{} \frac{\AAA}{(f)}\xrightarrow{} 0
\]
we get
\[
\dim_{\CC}\frac{\ann_{\AAA}(f)}{(c)} -\dim_{\CC}\frac{\AAA}{(f)}=\dim_{\CC}\frac{\AAA}{(c)}-\dim_{\CC}\AAA
\]
and the corollary follows.
\end{proof}

Now if $(X_1,\dots ,X_n)$ is a regular sequence our formula reduces to the formula of Gomez-Mont in the following way.
We find that
\[
\lambda =\dim_{\CC}\frac{\ann_{\BB}(f)}{(c)}=\dim_{\CC}\frac{\ann_{\BB}(c)}{(f)}
\]
and we get for even $n$
\begin{equation*}
\begin{split}
\IVn (X) &=\dim_{\CC}\frac{\BB}{(f)}-\dim_{\CC}\frac{\ann_{\AAA}(f)}{(c)}+\dim_{\CC}\frac{\AAA}{(c)}-\dim_{\CC}\AAA (c)\\
 &= \dim_{\CC}\frac{\BB}{(f)}-\dim_{\CC}\ann_{\BB}(f)\\
 &=\dim_{\CC}\frac{\BB}{(f)}-\dim_{\CC}\frac{\AAA}{(f)}.
\end{split}
\end{equation*}

For $n$ odd we get
\begin{equation*}
\begin{split}
\IVn (X) &=\dim_{\CC}\frac{\BB}{(f)}-\dim_{\CC}\frac{\ann_{\AAA}(f)}{(c)}+\dim_{\CC}\frac{\AAA}{(f)}\\
 &=\dim_{\CC}\frac{\BB}{(f)}+\dim_{\CC}\AAA (c)\\
 &=\dim_{\CC}\frac{\BB}{(f)}-\dim_{\CC}\frac{\AAA}{(c)}+\dim_{\CC}\AAA .
\end{split}
\end{equation*}
Note that the formula of Gomez-Mont can also be written in this form.

\section{Relations to residues and examples}
The index of the vector field $X$ can also be computed as a residue, see \cite{lss}:
Let $\delta_1 ,\dots ,\delta_{n-1}\in\RR_{>0}$ be chosen small enough and the real hypersurfaces $\{|X_i|=\delta_i\}$,
$i=1,\dots ,n-1$, in general position. Further let the real $(n-1)$-cycle
\[
\Sigma :=\{ f=0,|X_i|=\delta_i ,i=1,\dots ,n-1\}
\]
oriented so that $d(\arg X_1)\wedge\dots\wedge d(\arg X_{n-1})>0$. One defines the holomorphic function $\hat{c}$ to be
the coefficient of $t^{n-1}$ in the formal power series expansion of
\[
\frac{\det (\mathds{1}-\frac{ti}{2\pi} DX)}{\det (1-\frac{ti}{2\pi} c)}.
\]
Then one has
\[
\IVn (X)=\int_{\Sigma}\frac{\hat{c}dz_1\wedge\dots\wedge dz_{n-1}}{X_1\dots X_{n-1}}.
\]
$DX$ is the Jacobi matrix of $X$ here. Let us compute $\hat{c}$ explicitely.
Denote by $\sigma_i(DX)$ the coefficients of the characteristical polynomial of $DX$, i.e.
\[
\det (DX-tE)=\sum_{i=0}^n(-1)^i\sigma_{n-i}(DX)t^i.
\]
Set $g(t):=\det (E-\frac{ti}{2\pi} DX)$ and $h(t):=1/(1-\frac{ti}{2\pi} c)$. We have
\[
g^{(k)}(0)=\bigl(\frac{-i}{2\pi }\bigr)^kk!\sigma_k(DX)
\]
and
\[
h^{(l)}(0)=l!\bigl(\frac{ic}{2\pi }\bigr)^l.
\]
Then we get
\begin{equation*}
\begin{split}
\hat{c} &= \frac{(gh)^{n-1}(0)}{(n-1)!}\\
          &= \frac{1}{(n-1)!}\sum_{k=0}^{n-1}\binom{n-1}{k}g^{(k)}(0)h^{(n-k-1)}(0)\\
          &= \bigl(\frac{1}{2\pi i}\bigr)^{n-1}\sum_{k=0}^{n-1}(-1)^{n-k-1}c^{n-k-1}\sigma_k(DX)\\
          &= \bigl(\frac{1}{2\pi i}\bigr)^{n-1}\sum_{k=0}^{n-1}(-1)^kc^k\sigma_{n-k-1}(DX).\\
\end{split}
\end{equation*}
Using the formula in \cite{k} expressing integrals over cycles in $V$ as Grothendieck residues in $\CC^n$ we obtain
\[
\IVn (X)=\RCn
\begin{bmatrix}
\frac{\partial f}{\partial z_n}\sum_{k=0}^{n-1}(-1)^kc^k\sigma_{n-k-1}(DX)\\
X_1\dots X_{n-1} f
\end{bmatrix}.
\]
If moreover $(X_1,\dots ,X_n)$ is a regular $\OCn$-sequence we find that there is a $k\in\NN$ with $f^k\in I_n$. This means
that there is a $\beta\in\OCn$ with $f^k=\beta X_n\;\mo (I_{n-1})$ and therefore we have for any $h\in\OCn$
\begin{equation*}
\begin{split}
\RCn\begin{bmatrix}
hc\\ X_1\dots X_n
\end{bmatrix} &=\RCn\begin{bmatrix} hc\beta\\ X_1\dots X_{n-1} f^k\end{bmatrix}\\
\RCn\begin{bmatrix} h\frac{\partial f}{\partial z_n}\\ X_1\dots X_{n-1} f\end{bmatrix} &=
\RCn\begin{bmatrix} h\frac{\partial f}{\partial z_n}f^k\\ X_1\dots X_{n-1}f^{k+1}\end{bmatrix}\\
&= \RCn\begin{bmatrix} h\beta cf\\ X_1\dots X_{n-1}f^{k+1}\end{bmatrix}
\end{split}
\end{equation*}
which means that
\[
\RCn\begin{bmatrix}
hc\\ X_1\dots X_n
\end{bmatrix} =\RCn\begin{bmatrix} h\frac{\partial f}{\partial z_n}\\ X_1\dots X_{n-1} f\end{bmatrix}.
\]
We have obtained

\begin{proposition}
\label{comp}
(i) $\IVn (X)=\RCn
\begin{bmatrix}
\frac{\partial f}{\partial z_n}\sum_{k=0}^{n-1}(-1)^kc^k\sigma_{n-k-1}(DX)\\
X_1\dots X_{n-1} f
\end{bmatrix}$.\\
(ii) If $(X_1,\dots ,X_n)$ is a regular sequence
\[
\IVn (X)=\ICn (X)-\RCn\begin{bmatrix} \det (DX-c\mathds{1})\\ X_1\dots X_n\end{bmatrix}.
\]
\end{proposition}

We remark that the residues are very easy to compute: If $g_1,\dots g_n$ is a regular sequence and $h\in\OCn$ we can find
integers $k_1,\dots ,k_n$ with $(x_1^{k_1},\dots ,x_n^{k_n})\subset (g_1,\dots ,g_n)$ and therefore a matrix $A$ with
\[
\begin{pmatrix} x_1^{k_1}\\ \vdots \\ x_n^{k_n} \end{pmatrix} =
A\begin{pmatrix} g_1\\ \vdots \\ g_n \end{pmatrix}.
\]
If $d$ is the coefficient of $x_1^{k_1-1}\dots x_n^{k_n-1}$ in the power series expansion of $h\det A$ then
$\RCn \begin{bmatrix} h\\ g_1\dots g_n \end{bmatrix} = d$.

Let us consider an example.
$D_k\colon f=x^2y+y^{k-1}$ , $k\geq 4$.
\[
X:=\frac{k-2}{2(k-1)}x^{m+1}\frac{\partial}{\partial x}+\frac{1}{k-1}x^my\frac{\partial}{\partial y},\;\; m\geq 2.
\]
We have $c:=x^m$ here. The index can be computed easily with Proposition \ref{comp} and we obtain $\IVn (X)=(m-1)(k-1)$.
We want to verify that Corollary \ref{alg} gives the same value. Using residues again we compute
\[
\dim_{\CC}\frac{\OCn}{(X_1, f)}=(k-1)(m+1).
\]
Consider the monomial $x^my$ in $\frac{\OCn}{(x^{m+1}, x^2y+y^{k-1})}$. We find that the monomials $x^my,\dots ,x^my^{k-2}$ are
linearly independent in this algebra and therefore
\[
\dim_{\CC}\frac{\BCn}{(f)}=\dim_{\CC}\frac{\OCn}{(X_1, f)}-k+2 =m(k-1) +1.
\]
Now we claim
\[
\{ g\in\OCn :gf\in (X_1,X_2)\} =\OCn (x^m).
\]
Let
\[
gf=c_1x^{m+1} +c_2x^my.
\]
Then
\[
g(x^2y +y^{k-1})-c_2x^my\in\OCn (x^{m+1})
\]
and therefore
\[
g(x^2+y^{k-2})\in\OCn (x^m)
\]
since $y$ is not a zero divisor in $\frac{\OCn}{(x^{m+1})}$. This means $g\in\OCn (x^m)$ since
$x^2+y^{k-2}$ is not a zero divisor in $\frac{\OCn}{(x^m)}$. It follows that
\[
\dim_{\CC}\frac{\ann_{\BCn}(f)}{(c)}=0.
\]
We also find immediately $\dim_{\CC}\ACn =k$ and $\dim_{\CC}\frac{\ACn}{(c)}=k$. Since $\ann_{\BCnp}(X_2)=\BCnp (x)$ we find
\[
\frac{\ann_{\BCnp}(X_2)}{\ann_{\BCnp}(X_2)\cap\BCnp (f,\partial f/\partial y)}=
\frac{\BCnp (x)}{\BCnp (x)\cap \BCnp (f,\partial f/\partial y)}\cong\ACn' (x)
\]
with $\ACn'=\frac{\OCn}{ (f,\partial f/\partial y)}$. A basis of $\ACn'$ is given by
$1,x,xy,\dots ,xy^{k-2},y,\dots ,y^{k-2}$ and therefore $\dim_{\CC}\ACn' (x)=k$ since $y^{k-2}\sim x^2$ in $\ACn'$. Now
Corollary \ref{alg} gives in fact
\[
\IVn (X)=(m-1)(k-1).
\]

%The following environment contains the references

\end{document}